\documentclass[a4paper]{article}
\usepackage{graphicx}
\usepackage{amssymb}
\usepackage{amsmath}
\usepackage{amsfonts}
\usepackage{varioref}
\usepackage[bookmarksnumbered]{hyperref}
\usepackage[numbers]{natbib}

\title{ Penalty Robin--Robin domain decomposition schemes for contact problems of nonlinear elastic bodies}

\author{Ihor~I.~Prokopyshyn
\footnote{Pidstryhach Institute for Applied Problems of Mechanics and Mathematics,
Naukova 3-b, Lviv, 79060, Ukraine, ihor84@gmail.com, Corresponding author }\and
Ivan~I.~Dyyak \footnote{Ivan Franko National University of Lviv, Universytetska 1,
Lviv, 79000, Ukraine, dyyak@franko.lviv.ua}   \and
 Rostyslav~M.~Martynyak \footnote{Pidstryhach Institute for Applied Problems of Mechanics and Mathematics,
Naukova 3-b, Lviv, 79060, Ukraine, labmtd@iapmm.lviv.ua } \and Ivan~A.~Prokopyshyn
\footnote {Ivan Franko National University of Lviv, Universytetska 1, Lviv, 79000,
Ukraine, lviv.pi@gmail.com}}

\begin{document}

\maketitle

\begin{abstract}
In this paper we propose on continuous level several domain decomposition methods to
solve unilateral and ideal multibody contact problems of nonlinear elasticity. We also
present theorems about convergence of these methods.


\end{abstract}

\section{Introduction}
\label{prokopyshyn_contrib:1} Many domain decomposition techniques for contact problems
have been proposed on discrete level, particularly substructuring and FETI methods
\cite{Avery2009, Dostal2009}.

Domain decomposition methods (DDMs), presented in \cite{Bayada2004, Koko2003,
Krause2002, Sassi2008} for unilateral two-body contact problems of linear elasticity,
are obtained on continuous level. All of them require the solution of nonlinear
one-sided contact problems for one or both of the bodies in each iteration.

In works \cite{Prokopyshyn2008, Prokopyshyn2010, Dyyak2010b} we have proposed a class
of penalty parallel Robin--Robin domain decomposition schemes for unilateral multibody
contact problems of linear elasticity, which are based on penalty method and iterative
methods for nonlinear variational equations. In each iteration of these schemes we have
to solve in a parallel way some linear variational equations in subdomains.

In this contribution we generalize domain decomposition schemes, proposed in
\cite{Prokopyshyn2008, Prokopyshyn2010, Dyyak2010b} to the solution of unilateral and
ideal contact problems of nonlinear elastic bodies. We also present theorems about the
convergence of these schemes.

\section{Formulation of multibody contact problem} \label{prokopyshyn_contrib:2}

Consider a contact problem of $N$ nonlinear elastic bodies $\Omega _{\alpha } \subset
{\mathbb R}^{3} $ with sectionally smooth boundaries $\Gamma _{\alpha }$, $\alpha
=1,2,...,N$ (Fig.~1). Denote $\Omega =\bigcup _{\alpha =1}^{N}\Omega _{\alpha }  $.

\begin{figure}[h]
 \center{ \includegraphics[bb=0mm 0mm 208mm 296mm, width=63.156mm,
height=43.928mm,
 viewport=3mm 4mm 205mm 292mm]{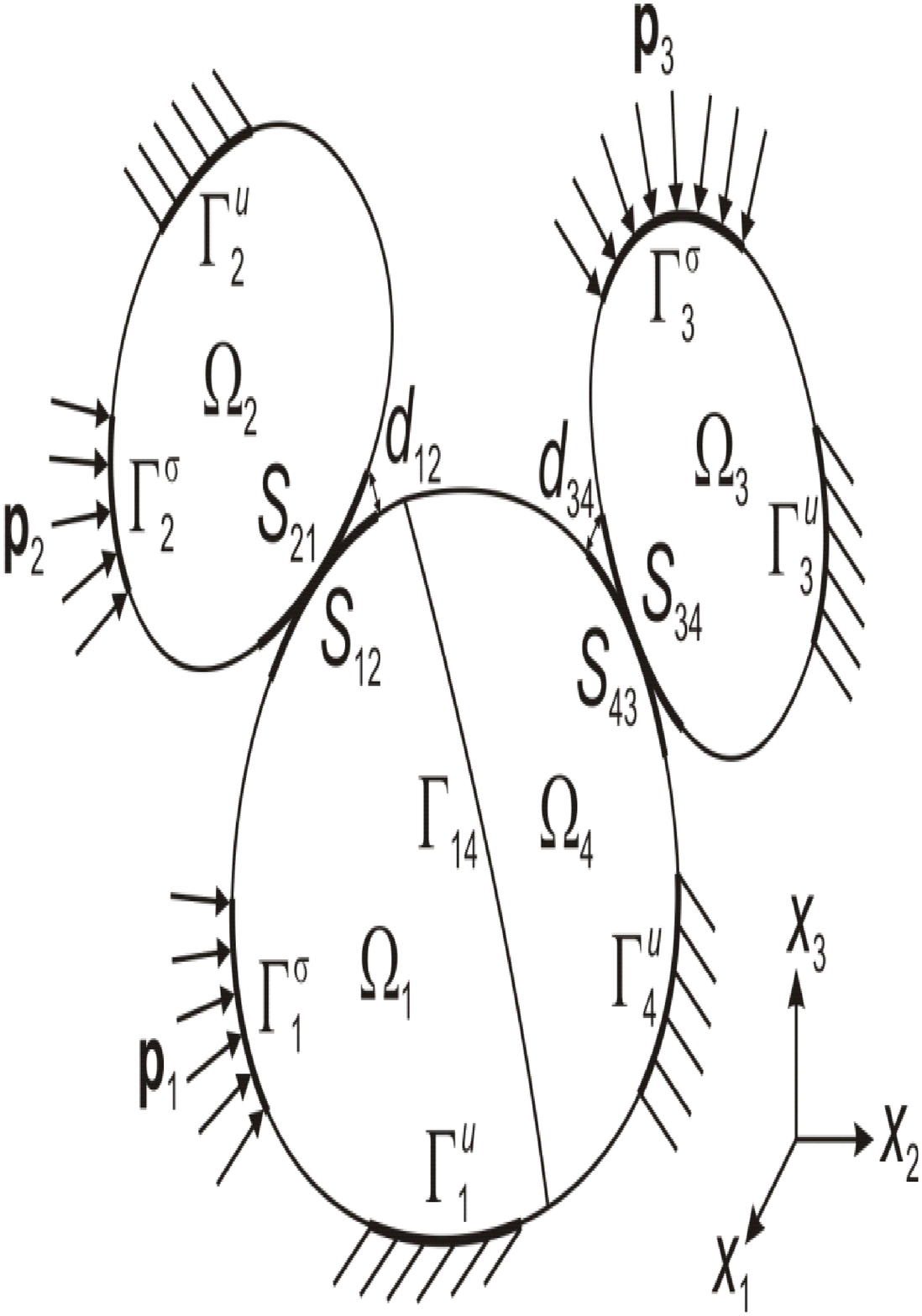}
 \\
 \textbf{Fig.~1} Contact of several bodies
 }
\end{figure}

A stress-strain state in point ${\bf x}=(x_{1} ,x_{2} ,x_{3})^{\top}$ of each body
$\Omega _{\alpha }$ is defined by the displacement vector
 ${\bf u}_{\,\alpha } =u_{\alpha \, i} \, {\bf e}_{i} $\,, the tensor of strains
 $ {\hat{\pmb{\varepsilon}}}_{\alpha }=\varepsilon _{\alpha \, ij}\,
 {\bf e}_{i} \, {\bf e}_{j} $
and the tensor of stresses
 $\hat{{\pmb \sigma }}_{\alpha } =\sigma _{\alpha \, ij} \, {\bf e}_{i} \,
 {\bf e}_{j} $\,.
These quantities satisfy Cauchy relations, equilibrium equations and nonlinear
 stress-strain law \cite{Ilyushin1948}:
\begin{equation} \label{prokopyshyn_contrib__1_}
 \sigma _{\alpha \, ij}   = \lambda _{\,\alpha }  \, \delta _{ij}
\,\Theta _{\alpha }  +2\, \mu _{\alpha }  \, \varepsilon _{\alpha \, ij}  -2\, \mu
_{\alpha }  \, \omega _{\alpha } (e_{\alpha } )\, e_{\alpha \, ij}\,, \, \, \,
i,j=1,2,3 \, \, ,
\end{equation}
\noindent where $ \Theta _{\alpha }  =\varepsilon _{\alpha \, 11} +\varepsilon _{\alpha
\, 22}  +\varepsilon _{\alpha \, 33} $ is the volume strain, $\lambda _{\,\alpha }({\bf
x})>0$,
 $\mu _{\alpha}({\bf x})>0$ are bounded Lame parameters,
 $e_{\alpha \, ij} =\varepsilon _{\alpha \,  ij} -{\delta _{ij} \,\Theta_{\alpha }
\mathord{\left/ {\vphantom {\delta _{ij} \Theta  3}} \right. \kern-\nulldelimiterspace}
3}$ are the components of the strain deviation tensor, $e_{\alpha} = {\sqrt{2 \,
g_{\alpha}}\, \mathord{\left/ {\vphantom {2\,g 3}}\right. \kern-\nulldelimiterspace}
3}$ is the deformation intensity, $ g_{\alpha } = (\varepsilon _{\alpha 11}
-\varepsilon _{\alpha 22} )^{2} +(\varepsilon _{\alpha 22} -\varepsilon _{\alpha 33}
)^{2} +(\varepsilon _{\alpha 33} -\varepsilon _{\alpha 11} )^{2} + 6\, (\varepsilon
_{\alpha 12}^{2} +\varepsilon _{\alpha 23}^{2} +\varepsilon _{\alpha 31}^{2})\,, $ and
$\omega _{\alpha } (z)$ is nonlinear differentiable function, which satisfies the
following properties:
\begin{equation} \label{prokopyshyn_contrib__2_}
0\le \omega _{\alpha } (z)\le {\partial \left(z\, \omega _{\alpha }
(z)\right)\mathord{\left/ {\vphantom {\partial \left(z\, \omega _{\alpha } (z)\right)
\partial z}} \right. \kern-\nulldelimiterspace} \partial z} <1\,, \,\,\,\, {\partial
\left(\omega _{\alpha } (z)\right)\mathord{\left/ {\vphantom {\partial \left(\omega
_{\alpha } (z)\right) \partial z}} \right. \kern-\nulldelimiterspace} \partial z} \ge 0
\,.
\end{equation}

On the boundary $\Gamma _{\alpha } $ let us introduce the local orthonormal basis
${\pmb \xi }_{\alpha } ,\, \, { \pmb \eta }_{\alpha } ,\, \, {\bf n}_{\,\alpha } $,
where $ {\bf n}_{\,\alpha } $ is the outer unit normal to $\Gamma _{\alpha } $. Then
the vectors of displacements and stresses on the boundary can be written in the
following way: ${\bf u}_{\,\alpha } =u_{\alpha \, \xi } \, { \pmb{\xi} }_{\alpha }
+u_{\alpha \eta } \, {\pmb{\eta} }_{\alpha } +u_{\alpha n} \, \bf{n}_{\,\alpha } , \,
\,
 {\pmb{ \sigma} }_{\alpha } =\hat{{\pmb{ \sigma }}}_{\alpha }
\cdot n_{\,\alpha } =\sigma _{\alpha \xi } \, {\pmb \xi }_{\alpha }
 +\sigma _{\alpha \eta } \, {\pmb \eta }_{\alpha } +\sigma _{\alpha \it n}
 \, n_{\,\alpha }\,.$

Suppose that the boundary $\Gamma _{\alpha } $ of each body consists of four disjoint
parts: $\Gamma _{\alpha } =\Gamma _{\alpha }^{u} \bigcup \Gamma _{\alpha }^{\sigma }
\bigcup \Gamma _{\alpha }^{I} \bigcup S_{\alpha } $, $\Gamma _{\alpha }^{u} \ne
\emptyset $, $\Gamma_{\alpha}^{u} = \overline{\Gamma_{\alpha}^{u}}$, $\Gamma _{\alpha
}^{I} \bigcup S_{\alpha } \ne \emptyset $, where $S_{\alpha } =\bigcup _{\beta \in
B_{\alpha } }S_{\alpha \beta }$, and $\Gamma _{\alpha }^{I} =\bigcup _{\beta ' \in
I_{\alpha } }\Gamma _{\alpha \beta '}$.
 Surface $S_{\alpha \beta } $ is the possible unilateral contact area of body
 $\Omega _{\alpha } $ with body $\Omega _{\beta } $,
 and $B_{\alpha } \subset \left\{1,2,...,N\right\}$
 is the set of the indices of all bodies in unilateral contact with body
 $\Omega _{\alpha } $. Surface $\Gamma _{\alpha \beta '} =\Gamma _{\beta '\alpha } $
 is the ideal contact area between bodies
  $\Omega _{\alpha } $ and $\Omega _{\beta '} $,
  and $I_{\alpha } \subset \left\{1,2,...,N\right\}$
  is the set of the indices of all bodies which have ideal contact with $\Omega _{\alpha } $.

We assume that the areas $S_{\alpha \beta } \subset \Gamma _{\alpha } $ and $S_{\beta
\alpha } \subset \Gamma _{\beta } $ are sufficiently close ($S_{\alpha \beta } \approx
S_{\beta \alpha } $), and  ${\bf n}_{\,\alpha } ({\bf x})\approx -{\bf n}_{\,\beta }
({\bf x}')$, ${\bf x}\in S_{\alpha \beta } $, ${\bf x}'=P({\bf x})\in S_{\beta \alpha
}$, where $P({\bf x})$ is the projection of ${\bf x}$ on $S_{\alpha \beta } $
\cite{Kravchuk1978}. Let $d_{\,\alpha \beta } ({\bf x})= {\pm \left\| {\bf x}-{\bf
x}'\right\|}_{2}$ be a distance between bodies $\Omega _{\alpha } $ and $\Omega _{\beta
} $ before the deformation. The sign of $d_{\,\alpha \beta }$ depends on  a statement
of the problem.

We consider homogenous Dirichlet boundary conditions on the part $\Gamma_{\alpha
}^{u}$, and Neumann boundary conditions on the part $\Gamma _{\alpha }^{\sigma}$:
\begin{equation} \label{prokopyshyn_contrib__3_}
{\bf u}_{\,\alpha } ({\bf x})=0, \,\,\, {\bf x}\in \Gamma _{\alpha}^{u}
\,;\,\,\,\,\,{\pmb \sigma }_{\alpha } ({\bf x})={\bf p}_{\alpha } ({\bf x}), \,\,\,
{\bf x}\in \Gamma _{\alpha }^{\sigma } \, .
\end{equation}

On the possible contact areas $S_{\alpha \beta } $, $\beta \in B_{\alpha}$, $\alpha
=1,2,...,N$ the following nonlinear unilateral contact conditions hold: \noindent
\begin{equation} \label{prokopyshyn_contrib__5_}
\sigma _{\alpha  n} ({\bf x})=\sigma _{\beta  n} ({\bf x}')\le 0 \,,\,\,\, \sigma
_{\alpha \, \xi } ({\bf x})=\sigma _{\beta \, \xi } ({\bf x}')=\sigma _{\alpha \, \eta
} ({\bf x})=\sigma _{\beta \, \eta } ({\bf x}')=0 \, ,
\end{equation}
\begin{equation} \label{prokopyshyn_contrib__7_}
u_{\alpha  n} ({\bf x})+u_{\beta  n} ({\bf x}')\le d_{\,\alpha \beta } ({\bf x}) \, ,
\end{equation}
\begin{equation} \label{prokopyshyn_contrib__8_}
\left(u_{\alpha  n} ({\bf x})+u_{\beta n} ({\bf x}')-d_{\,\alpha \beta } ({\bf x})\,
\right)\sigma _{\alpha  n} ({\bf x})=0 \, , \,\,\, {\bf x}\in S_{\alpha \beta } \, ,
\,\,\, {\bf x}'=P({\bf x}) \in S_{\beta \alpha } \, .
\end{equation}

On ideal contact areas $\Gamma _{\alpha \beta '} =\Gamma _{\beta ' \alpha } $, $\beta '
\in I_{\alpha}$, $\alpha =1,2,...,N$ we consider ideal mechanical contact conditions:
\begin{equation} \label{prokopyshyn_contrib__9_}
{\bf u}_{\,\alpha } ({\bf x})={\bf u}_{\,\beta '} ({\bf x})\,, \,\,\,{\pmb \sigma
}_{\alpha } ({\bf x})=-\, {\pmb \sigma }_{\beta '} ({\bf x}),\, \,\, {\bf x}\in \Gamma
_{\alpha \beta '} \, .
\end{equation}

\section{Penalty variational formulation of the \\ problem} \label{prokopyshyn_contrib:3}

For each body $\Omega _{\alpha } $ consider Sobolev space $V_{\alpha } =[H^{1} (\Omega
_{\alpha } )]^{3} $ and the closed subspace $V_{\alpha }^{0} =\left\{{\bf u}_{\,\alpha
} \in V_{\alpha }:\, \, {\bf u}_{\,\alpha } =0 \, \, \, {\rm on} \, \, \Gamma _{\alpha
}^{u} \right\}$. All values of the elements ${\bf u}_{\,\alpha } \in V_{\alpha } $,
${\bf u}_{\,\alpha } \in V_{\alpha }^{0} $ on the parts of boundary $\Gamma _{\alpha }$
should be understood as traces \cite{Kikuchi1988}.

Define Hilbert space $V_{0} =V_{1}^{0} \times ...\times V_{N}^{0} $ with the scalar
product $\left({\bf u}\,,{\bf v}\right)_{V_{0} } =\sum _{\alpha =1}^{N}\left({\bf
u}_{\,\alpha } ,{\bf v}_{\alpha } \right)_{V_{\alpha } }  $ and norm $\left\| {\bf
u}\right\| _{V_{0} } =\sqrt{\left({\bf u}\,,{\bf u}\right)_{V_{0} }}$, ${\bf u},{\bf
v}\in V_{0} $. Introduce the closed convex set of all displacements in $V_{0} $, which
satisfy nonpenentration contact conditions (\ref{prokopyshyn_contrib__7_}) and ideal
kinematic contact conditions:
\begin{equation}\label{prokopyshyn_contrib__11_}
K=\left\{{\bf u}\in V_{0}: \,  \, \,u_{\alpha \, n} +u_{\beta \, n} \le d_{\alpha \beta
}\, \, \,{\rm on}\, \, \,S_{\alpha \beta } ,\, \, \,{\bf u}_{\,\alpha '} ={\bf
u}_{\,\beta '} \, \, \, {\rm on}\, \, \, \Gamma _{\alpha '\beta '} \,\right\} ,
\end{equation}
where $\left\{\alpha ,\, \beta \right\}\in Q$, $Q=\left\{\left\{\alpha ,\beta
\right\}:\, \, \alpha \in \left\{1,2,...,N\right\},\,  \, \beta \in B_{\alpha }
\right\}$, $\left\{\alpha ',\, \beta '\right\}\in Q^{I}$, $Q^{I}=\left\{\left\{\alpha
',\beta '\right\}:\, \, \alpha '\in \left\{1,2,...,N\right\},\,  \, \beta '\in
I_{\alpha } \right\}$, and $d_{\alpha \beta } \in H_{\,00}^{1/2} ({\Xi}_{\alpha})$,
${\Xi}_{\alpha}={\rm {int}}\,({\Gamma}_{\alpha} \setminus {\Gamma}_{\alpha}^{u})$.

Let us introduce bilinear form $A({\bf u},{\bf v})=\sum _{\alpha=1}^{N}a_{\alpha} ({\bf
u}_{\,\alpha } ,{\bf v}_{\alpha } )$, ${\bf u},{\bf v}\in V_{0} $, which represents the
total elastic deformation energy of the system of bodies, linear form $L\,({\bf
v})=\sum _{\alpha =1}^{N}l_{\alpha } ({\bf v}_{\alpha } )$, ${\bf v}\in V_{0} $, which
is equal to the external forces work, and nonquadratic functional $H\,({\bf v})=\sum
_{\alpha =1}^{N}h_{\alpha } ({\bf v}_{\alpha })$, ${\bf v}\in V_{0} $, which represents
the total nonlinear deformation energy:
\begin{equation} \label{prokopyshyn_contrib__12_}
a_{\alpha } ({\bf u}_{\,\alpha } ,{\bf v}_{\alpha } )= \int _{\Omega _{\, \alpha } }
\,[\, \lambda _{\alpha } {\Theta_{\alpha }} ({\bf u}_{\,\alpha }) {\,\Theta_{\alpha }}
({\bf v}_{\alpha } )+2\, \mu _{\alpha } \sum _{i,j}\varepsilon _{\alpha \, ij} ({\bf
u}_{\,\alpha } )\, \varepsilon _{\alpha \, ij} ({\bf v}_{\alpha} ) \, \,]\,\, d\Omega
\, ,
\end{equation}
\begin{equation} \label{prokopyshyn_contrib__13_}
l_{\alpha } ({\bf v}_{\alpha } )=
 \int _{\Omega _{\alpha } }{\bf f}_{\,\alpha } \cdot {\bf v}_{\alpha }
 \,\,
 d\Omega  +\int _{\Gamma _{\alpha }^{\, \sigma } }
{\bf p}_{\alpha } \cdot {\bf v}_{\alpha } \,\, dS \, \, ,
\end{equation}
\begin{equation} \label{prokopyshyn_contrib__14_}
h_{\alpha } ({\bf v}_{\alpha } )=3\int _{\Omega _{\alpha } }\mu _{\alpha } \int _{\,
0}^{\, e_{\alpha } ({\bf v}_{\alpha } )}z\, \omega_{\alpha } (z)\,\, dz \,\,\, d\Omega
\, ,
\end{equation}
where ${\bf p}_{\alpha } \in [H_{\,00}^{-1/2} (\Xi_{\alpha})]^{3} $, and ${\bf
f}_{\,\alpha } \in [L_{2} (\Omega _{\alpha } )]^{3} $ is the vector of volume forces.

Using \cite{Kravchuk1978}, we have shown that the original contact problem has an
alternative weak formulation as the following minimization problem on the set $K$:
\begin{equation} \label{prokopyshyn_contrib__15_}
F({\bf u})=A\,({\bf u},{\bf u})/2-H({\bf u})-L\, ({\bf u})\to \mathop{\min
}\limits_{{\bf u}\, \in \, K} .
\end{equation}

Bilinear form $A$ is symmetric, continuous with constant $M_{A}>0$ and coercive with
constant $B_{A}>0$, and linear form $L$ is continuous. Nonquadratic functional $H$ is
doubly Gateaux differentiable in $V_{0} $:
\begin{equation} \label{prokopyshyn_contrib__4_}
H'({\bf u},{\bf v})=\sum _{\alpha}h'_{\alpha } ({\bf u}_{\,\alpha } ,{\bf v}_{\alpha }
),  \,\,\, H''({\bf u},{\bf v},{\bf w})=\sum _{\alpha}h''_{\alpha } ({\bf u}_{\,\alpha
} ,{\bf v}_{\alpha } ,{\bf w}_{\alpha } ) , \,\,\, {\bf u},{\bf v},{\bf w}\in V_{0} ,
\end{equation}
\begin{equation} \label{prokopyshyn_contrib__10_}
h'_{\alpha } ( {\bf u}_{\,\alpha } ,{\bf v}_{\alpha } )=2\int _{\Omega _{\alpha } }\,
\mu _{\alpha } \, \omega _{\alpha } (e_{\alpha } ({\bf u}_{\,\alpha } ))\, \sum
_{i,j}e_{\alpha  ij} ({\bf u}_{\,\alpha } )\,
 e_{\alpha  ij} ({\bf v}_{\alpha } ) \,\, d\Omega.
\end{equation}
Moreover, we have proved that the following conditions hold:
\begin{equation} \label{prokopyshyn_contrib__16_}
\left(\exists \,C>0\, \right)\left(\, \forall {\bf u}\in V_{0} \right)\,
 \left\{\,(1-C)\,A\,({\bf u},{\bf u}) \ge 2\,H\,({\bf u})\,\,\right\},
\end{equation}
\begin{equation} \label{prokopyshyn_contrib__17_}
\left(\forall {\bf u}\in V_{0} \right)\left(\exists R>0\, \right)\left(\forall {\bf
v}\in V_{0} \right)\left\{\, \left|H'({\bf u},{\bf v})\right|\le R\left\| {\bf
v}\right\| _{V_{0} } \right\},
\end{equation}
\begin{equation} \label{prokopyshyn_contrib__18_}
\left(\exists D>0\, \right)\left(\forall {\bf u},{\bf v},{\bf w}\in V_{0}
\right)\left\{\, \left|H''({\bf u},{\bf v},{\bf w})\right|\le D\left\| {\bf v}\right\|
_{V_{0} } \left\| {\bf w}\right\| _{V_{0} } \right\},
\end{equation}
\begin{equation} \label{prokopyshyn_contrib__19_}
\left(\exists B>0\, \right)\left(\forall {\bf u},{\bf v}\in V_{0} \right)\left\{A\,
({\bf v},{\bf v})-H''({\bf u},{\bf v},{\bf v})\ge B\left\| {\bf v}\right\| _{V_{0}
}^{2} \right\}.
\end{equation}

From these properties, it follows that there exists a unique solution $\bar{{\bf u}}
\in K$ of minimization problem (\ref{prokopyshyn_contrib__15_}), and this problem is
equivalent to the following variational inequality, which is nonlinear in ${\bf u}$:
\begin{equation} \label{prokopyshyn_contrib__20_}
A\, ({\bf u},{\bf v}-{\bf u})-H'({\bf u},{\bf v}-{\bf u})-L\, ({\bf v}-{\bf u})\ge
0,\,\,\, \forall \, {\bf v}\in K, \,\, {\bf u}\in K \, .
\end{equation}

To obtain a minimization problem in the whole space $V_{0} $, we apply a penalty method
\cite{Cea1971,Glowinski1976,Kikuchi1988,Lions1969} to problem
(\ref{prokopyshyn_contrib__15_}). We use a penalty in the form
\begin{eqnarray} \label {prokopyshyn_contrib__21_}
J_{\theta } ({\bf u})=\frac{1}{2\theta }  \sum _{\left\{\alpha ,\, \, \beta \right\}\,
\in \, Q}\left\| \left(d_{\alpha \beta } -u_{\alpha \, n} -u_{\beta \, n} \right)^{-}
\right\| _{L_{2} (S_{\alpha \beta } )}^{2} + \nonumber \\ + \frac{1}{2\theta } \sum
_{\left\{\alpha ',\, \, \beta '\right\}\, \in \, Q^{I} }\left\| {\bf u}_{\,\alpha '}
-{\bf u}_{\,\beta '} \right\| _{[L_{2} (\Gamma _{\alpha '\beta '} )]^{3} }^{2} ,
\end{eqnarray}
where $\theta >0$ is a penalty parameter, and $y^{-} =\min \{ 0,y\} $.

Now, consider the following unconstrained minimization problem in $V_{0} $:
\begin{equation} \label{prokopyshyn_contrib__22_}
F_{\theta } ({\bf u})=A\,({\bf u},{\bf u})/2-H\, ({\bf u})-L\, ({\bf u})+J_{\theta }
({\bf u})\to \mathop{\min }\limits_{{\bf u}\, \in \, V_{0} } .
\end{equation}

The penalty term $J_{\theta } $ is nonnegative and Gateaux differentiable in $V_{0} $,
and its differential $\,\,\,\,J'_{\theta } ({\bf u}, {\bf v})=-\frac{1}{\theta } \sum
_{\left\{\alpha ,\, \, \beta \right\}\, \in \, Q} \int _{S_{\alpha \beta}
}\left(d_{\alpha \beta }-u_{\alpha \, n}-u_{\beta \, n} \right)^{-}\, \left(v_{\alpha
\, n} +v_{\beta \, n} \right)\, dS +\frac{1}{\theta } \sum _{\left\{\alpha' ,\, \,
\beta' \right\}\, \in \, Q^{I}} \int _{{\Gamma}_{\alpha' \beta'} }\left({\bf
u}_{\alpha'}-{\bf u}_{\beta'} \right) \cdot \left({\bf v}_{\alpha'} - {\bf v}_{\beta'}
\right)\, dS$ satisfy the following properties \cite{Prokopyshyn2010}:
\begin{equation} \label{prokopyshyn_contrib__23_}
(\forall {\bf u}\in V_{0})(\exists \tilde{R}>0\,)(\forall {\bf v}\in V_{0})
\left\{\left|J'_{\theta }({\bf u},{\bf v})\right|\le \tilde{R}\left\| {\bf v}\right\|
_{V_{0} } \right\} ,
\end{equation}
\begin{equation} \label{prokopyshyn_contrib__24_}
(\exists \tilde{D}>0)(\forall {\bf u},{\bf v},{\bf w}\!\in\! V_{0})\! \left\{\!
\left|J'_{\theta } \, ({\bf u}+{\bf w},{\bf v})\!-\!J'_{\theta } \, ({\bf u},{\bf v})\,
\!\right|\!\le\! \tilde{D}\left\| {\bf v}\right\| _{V_{0}}\! \left\| {\bf w}\right\|
_{V_{0} }\!\right\},
\end{equation}
\begin{equation} \label{prokopyshyn_contrib__25_}
\left(\forall {\bf u},{\bf v}\in V_{0} \right)\left\{J'_{\theta } \,
 ({\bf u}+{\bf v},{\bf v})-J'_{\theta } \, ({\bf u},{\bf v})\ge 0\right\}.
\end{equation}

Using these properties and the results in \cite{Cea1971}, we have shown that problem
(\ref{prokopyshyn_contrib__22_}) has a unique solution $\bar{{\bf u}}_{\theta} \in
V_{0}$ and is equivalent to the following nonlinear variational equation in the space
$V_{0}$:
\begin{equation} \label{prokopyshyn_contrib__26_}
F'_{\theta } ({\bf u},{\bf v})=A\, ({\bf u},{\bf v})-H'({\bf u},{\bf v})+J'_{\theta }
({\bf u},{\bf v})-L\, ({\bf v})=0,\,\,\,
  \forall \, {\bf v}\in V_{0} , \,\,\, {\bf u}\in V_{0} .
\end{equation}

Using the results of works \cite{Glowinski1976, Lions1969}, we have proved that
$\left\| \bar{{\bf u}}_{\theta}-\bar{{\bf u}}\right\| _{V_{0} } \mathop{\to
}\limits_{\theta \to \,0 } 0$.

\section{Iterative methods for nonlinear variational \\ equations}

In arbitrary reflexive Banach space $V_{0}$ consider an abstract nonlinear variational
equation
\begin{equation} \label{prokopyshyn_contrib__27_}
\Phi \, ({\bf u}, {\bf v})=L\, ({\bf v}), \,\,\, \forall \, {\bf v}\in V_{0} , \,\,\,
{\bf u}\in V_{0}
\end{equation}
 where $\Phi :\, \, V_{0} \times
V_{0} \to {\mathbb R}$ is a functional, which is linear in ${\bf v}$, but nonlinear in
${\bf u}$, and $L$ is linear continuous form. Suppose that this variational equation
has a unique solution $\bar{{\bf u}}_{*} \in V_{0}$.

For the numerical solution of (\ref{prokopyshyn_contrib__27_}) we use the next
iterative method \cite{Dyyak2010a,Dyyak2010b,Prokopyshyn2010}:
\begin{equation} \label{prokopyshyn_contrib__28_}
G\, ({\bf u}^{k+1} ,{\bf v})=G\, ({\bf u}^{k} ,{\bf v})-\gamma ^{} \left[\Phi \, ({\bf
u}^{k} ,{\bf v})-L\, ({\bf v})\right], \,\,\,\forall \, {\bf v} \in V_{0}, \,\,\,\,
k=0,1,... \,\, ,
\end{equation}
where $G$ is some given bilinear form in $V_{0} \times V_{0} $, $\gamma \in {\mathbb
R}$ is fixed parameter, and ${\bf u}^{k} \in V_{0} $ is the \textit{k}-th approximation
to the exact solution of problem (\ref{prokopyshyn_contrib__27_}).

We have proved the next theorem \cite{Dyyak2010a,Prokopyshyn2010} about the convergence
of this method.

\textbf{Theorem~1.} \textit{Suppose that the following conditions hold}
\begin{equation} \label{prokopyshyn_contrib__29_}
\left(\forall {\bf u}\in V_{0} \right)\left(\exists R_{\Phi } >0\, \right)\left(\forall
{\bf v}\in V_{0} \right)\left\{\, \left|\Phi \, ({\bf u},{\bf v})\right|\le R_{\Phi }
\left\| {\bf v}\right\| _{V_{0} } \right\},
\end{equation}
\begin{equation} \label{prokopyshyn_contrib__30_}
\left(\exists D_{\Phi } \!>\!0\right)\!\left(\forall {\bf u},{\bf v},{\bf w}\!\in
\!V_{0} \right)\!\left\{\,\!\left|\Phi \, ({\bf u}+{\bf w},{\bf v})-\Phi \, ({\bf
u},{\bf v})\,\!\right| \le D_{\Phi } \!\left\| {\bf v}\right\| _{V_{0} } \!\left\| {\bf
w}\right\| _{V_{0} } \!\right\},
\end{equation}
\begin{equation} \label{prokopyshyn_contrib__31_}
\left(\exists B_{\Phi } >0\right)\left(\forall {\bf u},{\bf v}\in V_{0} \right)\,
\left\{\, \Phi \, ({\bf u}+{\bf v},{\bf v})-\Phi \, ({\bf u},{\bf v})\ge B_{\Phi }
\left\| {\bf v}\right\| _{V_{0} }^{2} \right\},
\end{equation}
\textit {bilinear form} $G$ \textit {is symmetric, continuous with constant} $M_{G}>0$
\textit {and coercive with constant} $B_{G}>0$, \textit {and}
 $\gamma \in (0;\,2\gamma^{*}) \, ,
  \gamma^{*} =B_{\Phi } {B_{G}} / {D_{\Phi}^{2}} $.

\textit{Then $\left\| {\bf u}^{k} -\bar{{\bf u}}_{*}\right\| _{V_{0} } \mathop{\to
}\limits_{k\to \infty } 0$, where $\{ {\bf u}^{k} \} \subset V_{0}$ is obtained by
method} (\ref{prokopyshyn_contrib__28_}). \textit {Moreover, the convergence rate in
norm $\left\| \, \cdot \, \right\| _{G} =\sqrt{G\, (\cdot ,\cdot )} $ is linear,}
\textit {and the highest convergence rate in this norm reaches as} $\gamma=\gamma^{*}$.

In addition, we have proposed nonstationary iterative method to solve
(\ref{prokopyshyn_contrib__27_}), where bilinear form $G$ and parameter $\gamma$ are
different in each iteration:
\begin{equation} \label{prokopyshyn_contrib__36_}
G^{k} ({\bf u}^{k+1} ,{\bf v})=G^{k} ({\bf u}^{k} ,{\bf v})-\gamma ^{k} \left[\Phi \,
({\bf u}^{k} ,{\bf v})-L\, ({\bf v})\right], \,\,\,\forall \, {\bf v} \in V_{0},
\,\,\,\, k=0,1,... \,\, .
\end{equation}
A convergence theorem for this method is proved in \cite{Prokopyshyn2010}.

\section{Domain decomposition schemes for contact \\ problems}

Now let us apply iterative methods (\ref{prokopyshyn_contrib__28_}) and
(\ref{prokopyshyn_contrib__36_}) to the solution of nonlinear penalty variational
equation (\ref{prokopyshyn_contrib__26_}) of multibody contact problem. This penalty
equation can be written in form (\ref{prokopyshyn_contrib__27_}), where
\begin{equation} \label{prokopyshyn_contrib__37_}
\Phi \, ({\bf u},{\bf v})=\, A\, ({\bf u},{\bf v})-H'({\bf u},{\bf v})+J'_{\theta }
({\bf u},{\bf v}), \,\,\,\, {\bf u},{\bf v}\in V_{0} .
\end{equation}

We consider such variants of methods (\ref{prokopyshyn_contrib__28_}) and
(\ref{prokopyshyn_contrib__36_}), which lead to the domain decomposition.

Let us take the bilinear form $G$ in iterative method (\ref{prokopyshyn_contrib__28_})
as follows \cite{Dyyak2010b,Prokopyshyn2010}:
\begin{equation} \label{prokopyshyn_contrib__38_}
G\, ({\bf u},{\bf v})=A\, ({\bf u},{\bf v})+X({\bf u},{\bf v}),\,\,\,\, {\bf u},{\bf
v}\in V_{0} ,
\end{equation}
$$
X({\bf u},{\bf v})=\frac{1}{\theta } \sum _{\alpha=1}^{N}\left[\sum _{\beta \in
B_{\alpha } }\int _{\, S_{\alpha \beta } }u_{\alpha \, n} v_{\alpha \, n} \,\psi
_{\alpha \beta } \, dS+\sum _{\beta ' \in I_{\alpha } }\int _{\, \Gamma _{\alpha \beta
'} }{\bf u}_{\alpha } \cdot {\bf v}_{\alpha } \,\phi _{\alpha \beta '} \, dS \right],
$$
where $\psi _{\alpha \beta } ({\bf x})\!=\!\{\,\!1,\, {\bf x}\in S_{\alpha \beta }^{1}
\,\}\vee \{\,\!0,\, {\bf x}\in S_{\alpha \beta } \backslash S_{\alpha \beta }^{1} \,\}$
and $\phi _{\alpha \beta '} ({\bf x})\!=\!\{\,\!1,\, {\bf x}\in \Gamma _{\alpha \beta
'}^{1} \,\}\vee \{\,\!0,\, \, \, {\bf x}\in \Gamma _{\alpha \beta '} \backslash \Gamma
_{\alpha \beta '}^{1} \,\}$ are characteristic functions of arbitrary subsets
$S_{\alpha \beta }^{1} \subseteq S_{\alpha \beta } $, $\Gamma _{\alpha \beta '}^{1}
\subseteq \Gamma _{\alpha \beta '} $ of possible unilateral and ideal contact areas
respectively.

Introduce a notation $\tilde{{\bf u}}^{k+1}=[{\bf u}^{k+1}-{\bf u}^{k}]/{\gamma}+{\bf
u}^{k} \in V_{0}$. Then iterative method (\ref{prokopyshyn_contrib__28_}) with bilinear
form (\ref{prokopyshyn_contrib__38_}) can be written in such way:
\begin{equation} \label{prokopyshyn_contrib__40_}
A \left(\tilde{{\bf u}}^{k+1} ,{\bf v}\right)+X \left(\tilde{{\bf u}}^{k+1} ,{\bf
v}\right)=L\, ({\bf v})+X\left({\bf u}^{k} ,{\bf v}\right)+H'({\bf u}^{k} ,{\bf
v})-J'_{\theta } ({\bf u}^{k} ,{\bf
v}), 
\end{equation}
\begin{equation} \label{prokopyshyn_contrib__41_}
{\bf u}^{k+1} =\gamma \,\, \tilde{{\bf u}}^{k+1} +\left(1-\gamma \right){\bf u}^{k} ,
\,\,\,\, k=0,1,...\,\,.
\end{equation}

Bilinear form $X$ is symmetric, continuous with constant $M_{X}>0$, and nonnegative
\cite{Prokopyshyn2010}. Due to these properties, and due to the properties of bilinear
form $A$, it follows that the conditions of Theorem~1 hold. Therefore, we obtain the
next proposition:

\textbf{Theorem~2.} \textit{The sequence} $\{{\bf u}^{k} \}$ \textit{of the method
(\ref{prokopyshyn_contrib__40_})~--~(\ref{prokopyshyn_contrib__41_}) converges strongly
to the solution of penalty variational equation} (\ref{prokopyshyn_contrib__26_})
\textit{for} $\gamma \in (0;\,2B_{\Phi } {B_{G}} / {D_{\Phi}^{2}})$, \textit{where}
$B_{G}=B_{A}$, $B_{\Phi}=B$, $D_{\Phi}=M_{A}+D+\tilde{D}$. \textit{The convergence rate
in norm $\left\| \, \cdot \, \right\| _{G}$ is linear}.

As the common quantities of the subdomains are known from the previous iteration,
variational equation (\ref{prokopyshyn_contrib__40_}) splits into $N$ separate
equations for each subdomain $\Omega _{\alpha } $, and method
(\ref{prokopyshyn_contrib__40_})~--~(\ref{prokopyshyn_contrib__41_}) can be written in
the following equivalent form:
\[a_{\alpha } (\tilde{{\bf u}}_{\,\alpha }^{k+1} ,
 {\bf v}_{\alpha } ) + \sum _{\beta
 \in B_{\alpha } } \int _{ S_{\alpha \beta } } \frac{\psi _{\alpha
 \beta}}{\theta }\, \tilde{u}_{\alpha\, n}^{k+1} \, v_{\alpha \, n}\,
 dS + \sum _{\beta' \in I_{\alpha }} \int _{
 \Gamma _{\alpha \beta'}}  \frac{\phi _{\alpha \beta '}}{\theta }\,
 \tilde{{\bf u}}_{\alpha }^{k+1} \cdot {\bf v}_{\alpha }
 \, dS  = \]
\[=l_{\alpha } ({\bf v}_{\alpha } )+\frac{1}{\theta } \sum _{\beta
\, \in B_{\alpha } }\int _{\, S_{\alpha \beta } }\left[\psi _{\alpha \beta} \,
u_{\alpha \, n}^{k}+\left(d_{\alpha \beta } -u_{\alpha \, n}^{k} -u_{\beta \, n}^{k}
\right)^{-}\right] v_{\alpha \, n} \, dS +  \]
\begin{equation} \label{prokopyshyn_contrib__42_}
+\frac{1}{\theta } \sum _{\beta ' \in I_{\alpha } }\int _{\, \Gamma _{\alpha \beta '}
}\left[\phi _{\alpha \beta '}\,{\bf u}_{\alpha }^{k}+\left({\bf u}_{\beta '}^{k} -{\bf
u}_{\alpha }^{k} \right)\right]\cdot {\bf v}_{\alpha } \, dS + h'_{\alpha } ( {\bf
u}_{\alpha}^{k} ,{\bf v}_{\alpha } ) \, , \,\,\,\,\forall \, {\bf v}_{\alpha} \in
V_{\alpha}^{0},
\end{equation}
\begin{equation} \label{prokopyshyn_contrib__43_}
{\bf u}_{\alpha }^{k+1} =\gamma \,\, \tilde{{\bf u}}_{\alpha }^{k+1} +\left(1-\gamma
\right){\bf u}_{\alpha }^{k}, \,\,\,\, \alpha =1,2,...,N, \,\,\,\, k=0,1, \, ...\,\,.
\end{equation}

In each iteration $k$ of method (\ref{prokopyshyn_contrib__42_}) --
(\ref{prokopyshyn_contrib__43_}), we have to solve $N$ linear variational equations in
parallel, which correspond to some linear elasticity problems in subdomains with
additional volume forces in $\Omega_{\alpha}\,$, and with Robin boundary conditions on
contact areas. Therefore, this method refers to parallel Robin--Robin type domain
decomposition schemes.

Taking different characteristic functions $\psi _{\alpha \beta }$ and $\phi _{\alpha
'\beta '}$, we can obtain different particular cases of penalty domain decomposition
method (\ref{prokopyshyn_contrib__42_})--(\ref{prokopyshyn_contrib__43_}).

Thus, taking $\psi _{\alpha \beta } ({\bf x})\equiv 0$, $\beta \in B_{\alpha } $, $\phi
_{\alpha \beta '} ({\bf x})\equiv 0$, $\beta '\in I_{\alpha } $, $\alpha =1,2,...,N$,
we get parallel Neumann--Neumann domain decomposition scheme.

Other borderline case is when $\psi _{\alpha \beta } ({\bf x})\equiv 1$, $\beta \in
B_{\alpha } $, $\phi _{\alpha \beta '}({\bf x})\equiv 1$, $\beta '\in I_{\alpha }$,
$\alpha =1,2,...,N$, i.e. $S_{\alpha \beta }^{1} =S_{\alpha \beta } $, $\Gamma _{\alpha
\beta '}^{1} =\Gamma _{\alpha \beta '} $.

Moreover, we can choose functions $\psi _{\alpha \beta } $ and $\phi _{\alpha \beta '}
$ differently in each iteration $k$. Then we obtain nonstationary domain decomposition
schemes, which are equivalent to iterative method (\ref{prokopyshyn_contrib__36_}) with
bilinear forms
\begin{equation} \label{prokopyshyn_contrib__44_}
G^{k}({\bf u},{\bf v})=A\, ({\bf u},{\bf v})+X^{k} ({\bf u},{\bf v}), \,\,\,\, {\bf
u},{\bf v}\in V_{0} , \,\,\,\, k=0,1,... \,\, ,
\end{equation}
\[X^{k}({\bf u},{\bf v})=\frac{1}{\theta } \sum _{\alpha=1}^{N}
\left[\sum _{\beta \in B_{\alpha } }\int _{\, S_{\alpha \beta } }u_{\alpha \, n}
v_{\alpha \, n} \,\psi _{\alpha \beta }^{k} \, dS+\sum _{\beta ' \in I_{\alpha } }\int
_{\, \Gamma _{\alpha \beta '} }{\bf u}_{\alpha } \cdot {\bf v}_{\alpha } \,\phi
_{\alpha \beta '}^{k} \, dS \right].\]

If we take characteristic functions $\psi _{\alpha \beta }^{k} $ and $\phi _{\alpha
\beta '}^{k} $ as follows \cite{Dyyak2010b,Prokopyshyn2008,Prokopyshyn2010}:
\[\psi _{\alpha \beta }^{k} ({\bf x})=\chi _{\alpha \beta }^{k}
({\bf x})=\left\{\begin{array}{c} {0,\, \, d_{\alpha \beta } ({\bf x})-u_{\alpha \,
n}^{k} ({\bf x})-u_{\beta \, n}^{k} ({\bf x}')\ge 0} \\ {1,\, \, \, d_{\alpha \beta }
 ({\bf x})-u_{\alpha \, n}^{k} ({\bf x})-u_{\beta \, n}^{k}
  ({\bf x}')<0} \end{array}\right. , \,\,{\bf x}'=P({\bf x}), \,\,{\bf x}\in S_{\alpha \beta } \,,
  \]
\[\phi _{\alpha \beta '}^{k} ({\bf x})\equiv 1,
\,\,\, {\bf x}\in \Gamma _{\alpha \beta '}, \,\,\, \beta \in B_{\alpha },\,\,\,\beta
'\in I_{\alpha },\,\,\,\alpha =1,2,...,N,\] then we shall get the method, which can be
conventionally named as nonstationary parallel Dirichlet--Dirichlet domain
decomposition scheme.

In addition to methods (\ref{prokopyshyn_contrib__28_}),
(\ref{prokopyshyn_contrib__38_}) and (\ref{prokopyshyn_contrib__36_}),
(\ref{prokopyshyn_contrib__44_}), we have proposed another family of DDMs for the
solution of (\ref{prokopyshyn_contrib__26_}), where the second derivative of functional
$H({\bf u})$ is used. These domain decomposition methods are obtained from
(\ref{prokopyshyn_contrib__36_}), if we choose bilinear forms $G^{k}({\bf u},{\bf v})$
as follows
\begin{equation} \label{prokopyshyn_contrib__48_}
G^{k}({\bf u},{\bf v})=A\, ({\bf u},{\bf v})-H''({\bf u}^{k},{\bf u},{\bf v})+X^{k}
({\bf u},{\bf v}), \,\,\,\, {\bf u},{\bf v}\in V_{0} , \,\,\,\, k=0,1,... \,\, .
\end{equation}

Numerical analysis of presented penalty Robin--Robin DDMs has been made for plane
unilateral two-body and three-body contact problems of linear elasticity
($\omega_{\alpha} \equiv 0$) using finite element approximations \cite{Dyyak2010b,
Prokopyshyn2008, Prokopyshyn2010}. Numerical experiments have confirmed the theoretical
results about the convergence of these methods.

Among the positive features of proposed domain decomposition schemes are the simplicity
of the algorithms and the regularization of original contact problem because of the use
of penalty method. These domain decomposition schemes have only one iteration loop,
which deals with domain decomposition, nonlinearity of the stress-strain relationship,
and nonlinearity of unilateral contact conditions.

\bibliography{prokopyshyn_contrib}

\end{document}